\documentclass[preprint,review,11pt]{elsarticle}

\usepackage{amssymb}
\usepackage{amsmath}
\usepackage{amsthm}
\usepackage{esint}
\newtheorem{theorem}{Theorem}[section]

\numberwithin{equation}{section}

\newtheorem{corollary}[theorem]{Corollary}

\newtheorem{lemma}[theorem]{Lemma}

\theoremstyle{definition}
\newtheorem{remark}[theorem]{Remark}

\journal{}

\usepackage{geometry}

\begin{document}

\baselineskip15pt

\begin{frontmatter}



\title{\huge Interior pointwise $C^{1,\alpha}$ estimates for Stokes systems in divergence form\tnoteref{t1}}
\tnotetext[t1]{Research supported by National Natural Science Foundation of China [grant 12001419 and grant 12071365] and the Fundamental Research Funds for the Central Universities [grant xzy012022004].}

\address[rvt1]{School of Mathematics and Statistics, Xi'an Jiaotong University, Xi'an 710049, China}

\author[rvt1]{Rong Dong}
\ead{dongrong1203@mail.xjtu.edu.cn }
\author[rvt1]{Dongsheng Li\corref{cor1}}
\ead{lidsh@mail.xjtu.edu.cn}
\cortext[cor1]{Corresponding author.}

\begin{abstract}
Interior pointwise $C^{1,\alpha}$ estimates are established for Stokes systems in divergence form where no continuity in time variable is assumed for the coefficients and the given data.
The estimates are attained by iteration and are presented by Campanato's characterization. The sharpness of the conclusions can be seen from Serrin's counterexample.
\end{abstract}

\bigskip


\end{frontmatter}


\section{Introduction}

The primary purpose of this paper is to establish interior pointwise $C^{1,\alpha}$ estimates for Stokes systems with variable coefficients in divergence form. Precisely, the system concerned here is the following
\begin{equation}\label{s0.0}
\left\{
\begin{array}{ll}
\vspace{2mm} u_t+\mathcal{L} u+\nabla p=\mbox{div}f~~\mbox{in}~~\Omega_{T};\\
\vspace{2mm}
\mbox{div}~u=g~~\mbox{in}~~\Omega_{T},
\end{array}
\right.
\end{equation}
 where
\begin{equation*}
\mathcal{L}:=-\mbox{div}\Big(A(x,t)\nabla \Big)=:-D_i\Big(a^{ij}(x,t)D_j\Big),
\end{equation*}
$\Omega_{T}=\Omega\times(-T,0]$ is a cylinder with $T>0$ and $\Omega\subset\mathbb{R}^{n}$ ($n\geq2$).
Throughout this paper,  the summation convention is assumed and the coefficient matrix $A(x,t)=\big(a^{ij}(x,t)\big)_{n\times n}$ always satisfies

\bigskip

($\bf H$) Ellipticity:
There exists $\lambda\in(0,1]$ such that
\begin{equation*}
a^{ij}(x,t)\xi_i\xi_j\geq\lambda|\xi|^2~~\mbox{and}~~|a^{ij}(x,t)|\leq\lambda^{-1}
\end{equation*}
for any $(x,t)\in\Omega_T$ and $\xi=(\xi_i)\in\mathbb{R}^{n}$.

\bigskip

For Stokes systems with constant coefficients, Solonnikov \cite{s11} established the global $C^{2,\alpha}$
(Schauder) estimates by using potential theory. Later, he further developed the global $C^{2,\alpha}$ estimates to Stokes systems in non-divergence form with variable coefficients (cf. \cite{s16} and the references therein).
In \cite{s1}, Abe and Giga showed the analyticity of the Stokes semigroup in $L^\infty_{\sigma}(\Omega)$, where the above estimates played an important part. For the Stokes systems in divergence form (1.1), the counterpart is the $C^{1,\alpha}$ estimates and this area seems remaining almost untouched. Here are some related interesting results with constant coefficients
and divergence free velocities.
Chang and Kang \cite{CK} solved the system in H\"{o}lder spaces by giving a global H\"{o}lder estimate of the velocities. Naumann and Wolf \cite{NW} showed that any weak solution is H\"{o}lder continuous with respect to the spatial variables that belong to $\mathbb{R}^3$.

Since Stokes systems are parabolic with degeneration, the global regularity heavily depends on the boundary data (cf. \cite{s1} \cite{CK} \cite{s11} and \cite{s16})
and the interior regularity is very different from the global ones.
This paper will investigate the interior $C^{1,\alpha}$ regularity for (1.1).
From Serrin's counterexample (cf. \cite{s14} or the following Remark 1.3 (ii)), one may find that the velocity is not H\"{o}lder continuous in the time variable even though the coefficients are constant and the righthand side terms are vanishing. That is, the regularity in the time direction is weak.
Actually, we establish the interior pointwise $C^{1,\alpha}$ regularity of velocities in the spatial variables provided $A$ and $f$ are only H\"{o}lder continuous in spatial variables, and $\nabla g$ is bounded. (This much improves the results in \cite{NW}.)
However an astonished result is that for curls of velocities, the H\"{o}lder continuity holds with respect to both spatial variables and time variable. (Observe no H\"{o}lder continuity in time variable is assumed for the coefficients and the righthand side terms.)
A corresponding result for $L_p$ estimates of (non-degenerate) parabolic equations is given by Krylov \cite{K}, where $W^{2}_p$ and $W^{1}_p$ estimates were presented for equations in non-divergence form and divergence form respectively
as the coefficients are only VMO in spatial variables. Such results (cf. Lemma 2.4 in the following) will be used in our proofs.
The conclusions in this paper are similar to directional regularity which one tries to reveal under conditions possessing regularity only on some directions.
Tian and Wang \cite{TW} obtained directional $C^{2,\alpha}$ regularity for elliptic equations, while Dong and Kim \cite{DK} further gave directional $C^{1,\alpha}$ and $C^{2,\alpha}$ estimates for both elliptic and parabolic equations. Dong, Li and Wang \cite{DLWd} and \cite{DLWj} studied directional homogenization problems and found different regularity between homogenization directions and non-homogenization directions.

Pointwise interior $C^{2,\alpha}$ and $C^{1,\alpha}$ estimates for fully nonlinear elliptic equations were originated by L.Caffarelli in his celebrated paper \cite{s8}.
Roughly speaking, under some structure regularity for the fully nonlinear elliptic operators, the previous estimate states that the solution will be $C^{2,\alpha}$ continuous at the same point if the nonlinear elliptic operator and the righthand side term are H\"{o}lder continuous at $x_0$; and the latter gives the $C^{1,\alpha}$ continuity of solutions as the oscillation of the nonlinear elliptic operator is of small BMO and the righthand side term is of $C^{\alpha-1}$ at the point concerned.
We will use the idea in \cite{s8} (cf. also \cite{s0}) to show the pointwise $C^{1,\alpha}$ estimate for the velocity of (1.1) in the space direction.  That is, the velocity is approximated at a $(1+\alpha)$-th order of the scales by a spatial variable first order polynomial whose coefficients may depend on the time. The approximation is indeed Campanato's characterization of $C^{1,\alpha}$ (cf. \cite{s10} \cite{GG} and \cite{L}) and the polynomial will be obtained by an iteration. In the iteration, the velocity, its gradient and its curl should be involved simultaneously and beyond expectation, for the curl of the velocity, the pointwise $C^{\alpha}$ estimate is achieved. We will estimate the curl of velocity by the gradient of it via parabolic equations, while estimate the gradient by the curl and velocity itself via elliptic equations.

The following are the main results of this paper.

\begin{theorem}
Suppose $M_0,M_1,M_2>0$, $0<\alpha<1$, $q>\frac{n}{1-\alpha}$,
Hypothesis $(H)$ holds,
and $u\in L_2(-T,0;W^{1}_2(\Omega))\cap L_\infty(-T,0;L_2(\Omega))$ and $p\in L_1(\Omega_T)$ satisfying (\ref{s0.0}) in weak sense.
If $Q_R(x_0,t_0)\subset\Omega_{T}$ and for any $0<r<R$,
\begin{equation}\label{t1}
\big|A(x,t)-A(x_0,t)\big|\leq M_0r^{\alpha}~~\mbox{in}~~Q_{r}(x_0,t_0),
\end{equation}
\begin{equation}\label{t2}
\fint_{Q_{r}(x_0,t_0)}|f-\overline{f}_{B_{r}(x_0)}|^2dxdt\leq M_1r^{2\alpha}
\end{equation}
and
\begin{equation}\label{t3}
\sup_{t\in(t_0-r^{2},t_0]}\fint_{B_{r}(x_0)}|\nabla g|^qdx\leq M_2,
\end{equation}
then there exist two constants $0<\sigma<1$ and $C>0$ depending only on $n$, $\lambda$, $\alpha$, $q$ and $M_0$, and two functions $a(t)\in L_\infty(t_0-\sigma^2 R^2,t_0;\mathbb{R}^{n})$ and $b(t)\in L_\infty(t_0-\sigma^2 R^2,t_0;\mathbb{R}^{n^2})$ such that
\begin{equation*}
\begin{aligned}
||a||_{L_\infty(t_0-\sigma^2 R^2,t_0)}+
||b||_{L_\infty(t_0-\sigma^2 R^2,t_0)}\leq
 C{\cal K},
 \end{aligned}
\end{equation*}
\begin{equation}\label{p01}
\begin{aligned}
\fint_{Q_{r}(x_0,t_0)}|\nabla\times u-\overline{\nabla\times u}_{Q_{r}(x_0,t_0)}|^2dxdt
\leq Cr^{2\alpha}{\cal K}^2
\end{aligned}
\end{equation}
and
\begin{equation}\label{p02}
\begin{aligned}
\sup_{t\in(t_0-r^{2},t_0]}\fint_{B_{r}(x_0)}|u-a(t)-b(t)x|^2dx
\leq&\ Cr^{2(1+\alpha)}{\cal K}^2
\end{aligned}
\end{equation}
for any $0<r\leq\sigma R$, where
$${\cal K}=\|u\|_{L_\infty(t_0-R^2,t_0;L_2(B_R(x_0)))}+\|\nabla\times u\|_{L_2(Q_R(x_0,t_0))}+M^\frac{1}{2}_1+M^\frac{1}{q}_2.$$
\end{theorem}

\bigskip

From the above pointwise results, the classical $C^{1,\alpha}$ estimates follow directly.

\bigskip

\begin{corollary}
Suppose $M_0>0$, $0<\alpha<1$,
and $A\in C_{x}^{\alpha}(\Omega_{T})$ with
$[A]_{x,\alpha;\Omega_T}\leq M_0$ and that $(H)$ holds.
Let
$u\in L_2(-T,0;W^{1}_2(\Omega))\cap L_\infty(-T,0;L_2(\Omega))$, $p\in L_1(\Omega_T)$, $f\in C_{x}^{\alpha}(\Omega_{T};\mathbb{R}^{n^2})$ and $g\in L_\infty(-T,0;W^{1}_\infty(\Omega))$ satisfying (\ref{s0.0}) in weak sense. Then for any
$0<T'<T$ and any $\Omega'\subset\subset\Omega$, we have
$\nabla\times u\in C^{\alpha,\frac{\alpha}{2}}(\overline{\Omega'_{T'}})$ and
$\nabla u\in C_x^{\alpha}(\overline{\Omega'_{T'}})$ such that
\begin{equation}\label{g01}
\begin{aligned}
&\|\nabla\times u\|_{C^{\alpha,\frac{\alpha}{2}}(\overline{\Omega'_{T'}})}+\|\nabla u\|_{C_x^{\alpha}(\overline{\Omega'_{T'}})}\\
&\mbox{}\hskip0.5cm\leq
C\Big(\|u\|_{L_\infty(-T,0;L_2(\Omega))}+\|\nabla\times u\|_{L_2(\Omega_{T})}
+\|f\|_{C_{x}^{\alpha}(\Omega_{T})}
+\|\nabla g\|_{L_\infty(\Omega_T)}\Big)
\end{aligned}
\end{equation}
where $C$ depends only on $n$, $\lambda$, $\alpha$, $M_0$ and $\min\{\mbox{dist}(\Omega',\partial\Omega),\sqrt{T-T'}\}$.

\end{corollary}

\bigskip

\begin{remark}
(i) Theorem 1.1 and Corollary 1.2 are effectively $C^{1,\alpha}$ regularity for weak solutions of Stokes systems in divergence form.

(ii) Let $h$ be a harmonic function in $\Omega$ and $\sigma$ be a function of $L_1(-T,0)$. Then $u(x,t)=\nabla h(x)\int_{-T}^t\sigma(s)ds$ and $p(x,t)=-h(x)\sigma(t)$ is a pair of weak solution of (1.1) with vanishing $f$ and $g$.
(This counterexample was given by Serrin \cite{s14}.)
It is easy to see that $Du=D\nabla h(x)\int_{-T}^t\sigma(s)ds$, which may not be H\"{o}lder continuous in variable $t$. That is, our conclusions are sharp.

(iii) Observe that $\nabla\times u$ has H\"{o}lder continuity in $t$ although $A$, $f$ and $g$ are not assumed to be continuous in $t$.


\end{remark}

\bigskip

The remainder of the paper will be organized as follows. In Section 2, we give some preliminary tools. The proof of Theorem 1.1 and Corollary 1.2 will be presented in Section 3.

\bigskip

Before the end of this section, we introduce some notations.
\vskip 2mm\noindent{\bf Notations:}

$B_r(x):=\{y\in\mathbb{R}^n:|x-y|<r\}$;
~$B_r:=B_r(0)$.
\vspace{2mm}

$Q_r(x,t):=B_r(x)\times(t-r^2,t]$;~
$Q_r:=B_r\times(-r^2,0]$.
\vspace{2mm}

$\overline{u}_{Q_r(x,t)}:=\fint_{Q_r(x,t)}u(y,s)
dyds=\frac{1}{|Q_r(x,t)|}\int_{Q_r(x,t)}u(y,s)
dyds$.
\vspace{2mm}

$\overline{u}_{B_r(x)}(s):=\fint_{B_r(x)}u(y,s)
dy=\frac{1}{|B_r(x)|}\int_{B_r(x)}u(y,s)dy$.
\vspace{2mm}

$D_iu:=\frac{\partial u}{\partial x_i}$; $Du:=\big(D_{1}u,D_{2}u,...,D_{n}u\big)$; $D^2u:=\big(D_iD_ju\big)_{n\times n}$.
\vspace{2mm}

$\nabla\times u:=\big(D_ju_i-D_iu_j\big)_{n\times n}$.
\vspace{2mm}

$\partial_p \Omega_T: \mbox{parabolic~boundary~of}~\Omega_T$.
\vspace{2mm}

$[u]_{\alpha,\frac{\alpha}{2};\Omega_T}:=\sup_{(x,t),(y,s)\in \Omega_T, (x,t)\neq (y,s)}\frac{|u(x,t)-u(y,s)|}{|x-y|^{\alpha}+|t-s|^{\frac{\alpha}{2}}}$.
\vspace{2mm}

$\|u\|_{C^{\alpha,\frac{\alpha}{2}}(\Omega_T)}:=[u]_{\alpha,\frac{\alpha}{2};\Omega_T}+\|u\|_{L_\infty(\Omega_T)}$.

\vspace{2mm}

$C^{\alpha,\frac{\alpha}{2}}(\Omega_T):=\{u:\|u\|_{C^{\alpha,\frac{\alpha}{2}}(\Omega_T)}<\infty\}$.
\vspace{2mm}

$[u]_{x,\alpha;\Omega_T}:=\sup_{(x,t),(y,t)\in \Omega_T, x\neq y}\frac{|u(x,t)-u(y,t)|}{|x-y|^{\alpha}}$.
\vspace{2mm}

$\|u\|_{C_{x}^\alpha(\Omega_T)}:=[u]_{x,\alpha;\Omega_T}+\|u\|_{L_\infty(\Omega_T)}$.

\vspace{2mm}

$C^\alpha_{x}(\Omega_T):=\{u:\|u\|_{C_{x}^\alpha(\Omega_T)}<\infty\}$.

\section{Preliminary}
We first introduce the lemma below, which characterizes the H\"{o}lder continuity by Campanato spaces.

\bigskip

\begin{lemma}
Let $T>0$, $0<\alpha<1$ and $\Omega\subset\mathbb{R}^n$ be a bounded Lipschitz domain.

(i) If
\begin{eqnarray*}
\sup_{(x,t)\in\Omega_{T}}\sup_{r>0}
\Big(\fint_{Q_r(x,t)\cap\Omega_{T}}|u(y,s)-\overline{u}_{Q_r(x,t)\cap\Omega_{T}}|^2dyds\Big)^{\frac{1}{2}}\leq
Ar^{\alpha}
\end{eqnarray*}
for some constant $A$,
then $u\in C^{\alpha,\frac{\alpha}{2}}(\overline{\Omega_{T}})$ and
$$[u]_{\alpha,\frac{\alpha}{2};\overline{\Omega_{T}}}\leq C
A$$
for some constant $C$ depending only on $\alpha$, $n$ and $\Omega_T$.

(ii) If there exist constants $r_0, A>0$ and functions $a(x,t)$ and $b(x,t)$ bounded in $\Omega_T$ such that
\begin{eqnarray*}
\sup_{0<r\leq r_0}
\Big(\fint_{B_r(x)\cap\Omega}|u(y,t)-a(x,t)-b(x,t)y|^2dy\Big)^{\frac{1}{2}}\leq
Ar^{1+\alpha}
\end{eqnarray*}
for any $(x,t)\in\Omega_T$,
then $\nabla u\in C_x^{\alpha}(\overline{\Omega_{T}})$ and
$$[\nabla u]_{x,\alpha;\overline{\Omega_{T}}}\leq C
A$$
for some constant $C$ depending only on $\alpha$, $n$, $T$ and $\Omega$.
\end{lemma}

\bigskip

We refer to \cite{L} for the proof of Lemma 2.1 (i), which is for anisotropic spaces, and to \cite{s10} for the proof of Lemma 2.1 (ii). The next Lemma is actually an equivalence between norms of polynomials, which may be well known.

\bigskip

\begin{lemma}
Let $T>0$, $a(t)\in L_{\infty}(-T,0;\mathbb{R})$ and $b(t)\in L_{\infty}(-T,0;\mathbb{R}^n)$.
If there exist $A$ and $r>0$ such that
\begin{equation}\label{j06}
\begin{aligned}
\sup_{t\in(-T,0]}\left(\fint_{B_{r}}|a(t)-b(t)x|^2dx\right)^{\frac12}
\leq A,
\end{aligned}
\end{equation}
then
\begin{equation}\label{j07}
\begin{aligned}
\sup_{t\in(-T,0]}|a(t)|
\leq CA~~\mbox{and}~~\sup_{t\in(-T,0]}|b(t)|
\leq CA/r
\end{aligned}
\end{equation}
for some constant $C$ depending only on $n$.
\end{lemma}

\vskip 2mm\noindent{\bf Proof.}\quad
Denote $P_1$ the space of one order polynomials and suppose $a+bx\in P_1$ with $a\in\mathbb{R}$ and $b,x\in\mathbb{R}^n$. Since both $\max\{|a|,|b|\}$ and $\left(\int_{B_1}|a+bx|^2dx\right)^{\frac12}$ are norms on $P_1$,
we have they are equivalent, that is, there exists a constant $C$ depending only on $n$ such that
$$C^{-1}\left(\int_{B_1}|a+bx|^2dx\right)^{\frac12}\leq\max\{|a|,|b|\}\leq C\left(\int_{B_1}|a+bx|^2dx\right)^{\frac12}.$$
It follows that
$$C^{-1}\left(\fint_{B_r}|a+bx|^2dx\right)^{\frac12}\leq\max\{|a|,r|b|\}\leq C\left(\fint_{B_r}|a+bx|^2dx\right)^{\frac12}$$
for any $r>0$. Then (\ref{j07}) is a consequence of (\ref{j06}).
\qed

\bigskip

The last two lemmas are $L_p$ estimates for elliptic and parabolic equations and we include them here for convenience.

\bigskip

\begin{lemma}
Let $1<p<\infty$ and $\tilde\Omega\subset\subset\Omega$.

(i) If $u\in W^1_{p}(\Omega;\mathbb{R}^{n})$, then for any $\zeta\in\mathbb{R}^{n}$ and $\eta\in \mathbb{R}$ independent of $x$,
\begin{eqnarray*}
\|\nabla u\|_{L_p(\tilde{\Omega})}\leq C\Big(\|\nabla\times u-\zeta\|_{L_p(\Omega)}
+\|\mbox{div}~u-\eta\|_{L_p(\Omega)}+\|u\|_{L_1(\Omega)}\Big),
\end{eqnarray*}
where $C$ depends only on $n$, $p$, $\tilde\Omega$ and $\Omega$.

(ii) If $u\in W^2_{p}(\Omega;\mathbb{R}^{n})$, then
\begin{eqnarray*}
\|D^2u\|_{L_p(\tilde{\Omega})}\leq C\Big(\|\mbox{div}~(\nabla\times u)\|_{L_p(\Omega)}
+\|\nabla(\mbox{div}~u)\|_{L_p(\Omega)}+\|u\|_{L_1(\Omega)}\Big),
\end{eqnarray*}
where $C$ depends only on $n$, $p$, $\tilde\Omega$ and $\Omega$.
\end{lemma}
\vskip 2mm\noindent{\bf Proof.}\quad
Recall $\nabla\times u=D_ju_i-D_iu_j$ and we have the following identity
\begin{eqnarray*}
\begin{aligned}
\Delta u
=\mbox{div}~\nabla\times u+\nabla~\mbox{div}~u
=\mbox{div}~(\nabla\times u-\xi)+\nabla~(\mbox{div}~u-\eta)
~~\mbox{in}~~\Omega
\end{aligned}
\end{eqnarray*}
for any $\zeta\in\mathbb{R}^{n}$ and $\eta\in \mathbb{R}$ independent of $x$.
Then by $W^1_p$ and $W^2_p$ estimates for elliptic equations, we obtain (i) and (ii) immediately.
\qed

\bigskip

\begin{lemma}
Let $2\leq p<\infty$, $0<T'<T$ and $\tilde\Omega\subset\subset\Omega\subset\mathbb{R}^n$ be bounded domains. Suppose  $(a^{ij}(t))_{n\times n}$ (not depending on $x$) satisfies $(H)$. Consider the equation:
\begin{equation}\label{j10}
u_t-a^{ij}(t)D_{ij}u=g~~\mbox{in}~~\Omega_T.
\end{equation}

(i) If $g=0$ and $u$ is a weak solution of (\ref{j10}), then $u, Du, u_t\in L_{\infty,loc}(\Omega_T)$ and
\begin{eqnarray*}
\|u\|_{L_\infty(\tilde{\Omega}_{T'})}+
\|Du\|_{L_\infty(\tilde{\Omega}_{T'})}+\|u_t\|_{L_\infty(\tilde{\Omega}_{T'})}\leq C\|u\|_{L_p(\Omega_{T})},
\end{eqnarray*}
where $C$ depends only on $n$, $\lambda$, $p$, $T'$, $T$, $\tilde\Omega$ and $\Omega$.

(ii) Assume that $\Omega$ is smooth. If $g=D_{ij}f_{ij}$ for some $f=(f_{ij})\in L_{2}(\Omega_{T};\mathbb{R}^{n^2})$ and $u\in L_{2}(\Omega_{T})$ is the generalized solution of (\ref{j10}) with
$u=0~\mbox{on}~\partial_p\Omega_{T}$ (See the following Remark 2.5), then
\begin{eqnarray}\label{j11}
\|u\|_{L_2({\Omega}_T)}\leq C\|f\|_{L_2(\Omega_{T})},
\end{eqnarray}
where $C$ depends only on $n$, $\lambda$ and $\Omega_{T}$.
\end{lemma}

\vskip 2mm\noindent{\bf Proof.}\quad
Let $\Omega^*\subset\Omega$ with $dist(\Omega^*,\partial \Omega)>\varepsilon_0>0$
 and $\tilde\Omega\subset\subset\Omega^*$.
 For $0<\varepsilon<\varepsilon_0$, denote
$\rho_\varepsilon$ the scaled mollifier in $\mathbb{R}^{n}$.
Then in view of $(a^{ij}(t))$ not depending on $x$, $u_\varepsilon=\rho_\varepsilon\ast u$ satisfies
\begin{equation}\label{r.15}
(u_{\varepsilon})_t-a^{ij}(t)D_{ij}u_\varepsilon=0~~\mbox{in}~~\Omega^*_{T}.
\end{equation}
Differentiating the equation (\ref{r.15}) suitable times with respect to $x$, and by the interior $L_p$ estimates for parabolic equations in \cite{K} and embedding inequalities, we have
\begin{eqnarray*}
\|u_{\varepsilon}\|_{L_\infty(\tilde{\Omega}_{T'})}+\|u_{\varepsilon t}\|_{L_\infty(\tilde{\Omega}_{T'})}+
\|Du_\varepsilon\|_{L_\infty(\tilde{\Omega}_{T'})}
\leq C\|u_\varepsilon\|_{L_p(\Omega^*_{T})},
\end{eqnarray*}
where $C$ depends only on $n$, $\lambda$, $p$, $\tilde\Omega_{T}$ and $\Omega^*_{T}$.
Then (i) is obtained by letting $\varepsilon\rightarrow0$.

Now we turn to the proof of  (ii). Suppose $\{f_{k}\}\subset W^{2,1}_2(\Omega_T)$ such that $f_k\rightarrow f$ in $L_2$.
Then (\ref{j10}) has unique solution sequence $u_k\in W^{2,1}_2(\Omega_T)$ with $u_k=0~\mbox{on}~\partial_p\Omega_{T}$. By $L_2$ estimates for parabolic equation, (\ref{j11}) holds for $u_k$. Let $k\rightarrow\infty$ and we have (\ref{j11}).
\qed

\bigskip

\begin{remark}
In the proof of (ii), it is easy to see the limit $u$ does not depend on the choice of $\{f_{k}\}$. We call this $u$ the generalized solution of (\ref{j10}) with $u=0~\mbox{on}~\partial_p\Omega_{T}$.
\end{remark}

\bigskip

\section{Proofs of main results}\setcounter{equation}{0}
To prove Theorem 1.1, we need only to show:

\begin{theorem}
Let $T>0$, $0<\alpha<1$ and  $q>\frac{n}{1-\alpha}$. Suppose $Q_1\subset\Omega_{T}$, $(H)$ holds, and $u\in L_2(-T,0;W^{1}_2(\Omega))\cap L_\infty(-T,0;L_2(\Omega))$ and $p\in L_1(\Omega_T)$ satisfying (\ref{s0.0}) in weak sense. There exist constants $0<r_0<1$, $0<\varepsilon_0<1$ and $C>0$ depending only on $n$, $\lambda$, $\alpha$ and $q$ such that if for any $0<r\leq1$,
\begin{equation}\label{sh.23}
\big|A(x,t)-A(0,t)\big|\leq\varepsilon_0r^{\alpha}~~\mbox{in}~~Q_r,
\end{equation}
\begin{equation}\label{sh.24}
\sup_{t\in(-1,0]}\fint_{B_1}|u|^2dx\leq1,~~\fint_{Q_1}|\nabla\times u|^2dxdt\leq1,
\end{equation}
\begin{equation}\label{n.25}
\begin{aligned}
\fint_{Q_r}|f-\overline{f}_{B_r}|^2dxdt\leq\varepsilon^2_0r^{2\alpha}~~~\mbox{and}
\end{aligned}
\end{equation}
\begin{equation}\label{sh.25}
\begin{aligned}
\sup_{t\in(-r^2,0]}\fint_{B_r}|\nabla g|^qdx\leq1,
\end{aligned}
\end{equation}
then there are $a(t)\in L_\infty(-r^2_0,0;\mathbb{R}^{n})$ and $b(t)\in L_\infty(-r^2_0,0;\mathbb{R}^{n^2})$ satisfying
\begin{equation}\label{j01}
\sup_{t\in(-r_0^{2},0]}|a(t)|~~\mbox{and}~~\sup_{t\in(-r_0^{2},0]}|b(t)|\leq C,
\end{equation}
\begin{equation}\label{sh.76}
\fint_{Q_{r}}|\nabla\times u-\overline{(\nabla\times u)}_{Q_{r}}|^2dxdt\leq Cr^{2\alpha}~~\mbox{and}
\end{equation}
\begin{equation}\label{f4.09}
\begin{aligned}
\sup_{t\in(-r^{2},0]}\fint_{B_{r}}|u-a(t)-b(t)x|^2dx
\leq Cr^{2(1+\alpha)}
\end{aligned}
\end{equation}
for any $0<r\leq r_0$.
\end{theorem}

Actually, Theorem 1.1 follows from Theorem 3.1 by the following way of normalization. After suitable choice of coordinates, we may suppose $(x_0,t_0)=(0,0)$.

Choose $0<\tau<R$ such that $\tau^\alpha M_0\leq\varepsilon_0$.
For $(x,t)\in Q_1$ and $0<r\leq1$, set
$$\tilde{A}(x,t)=A(\tau x,\tau^2t),~~\tilde{u}(x,t)=\frac{u(\tau x,\tau^2t)}{\tau J},~~\tilde{p}(x,t)=\frac{p(\tau x,\tau^2t)}{J},$$
$$\tilde{f}(x,t)=\frac{f(\tau x,\tau^2t)}{J}~~\mbox{and}~~\tilde{g}(x,t)=\frac{g(\tau x,\tau^2t)}{J},$$
where
\begin{eqnarray*}
\begin{aligned}
J=&\ \frac{1}{\tau|B_\tau|^{\frac12}}\|u\|_{L_\infty(-T,0;L_2(\Omega))}+\frac{1}{|Q_\tau|^{\frac12}}\|\nabla\times u\|_{L_2(\Omega_T)}+\frac{1}{\varepsilon_0}M^{\frac12}_1+M^\frac{1}{q}_2.
\end{aligned}
\end{eqnarray*}
Then $\tilde u$, $\tilde p$, $\tilde f$ and $\tilde g$ satisfy (1.1), 
\begin{equation*}
\big|\tilde{A}(x,t)-\tilde{A}(0,t)\big|=\big|A(\tau x,\tau^2t)-A(0,\tau^2t)\big|\leq M_0\tau^{\alpha}r^{\alpha}\leq\varepsilon_0r^{\alpha}~~\mbox{in}~~Q_r,
\end{equation*}
\begin{equation*}
\sup_{t\in(-1,0]}\fint_{B_1}|\tilde{u}|^2dx
=\frac{1}{\tau^2J^2}\sup_{t\in(-\tau^2,0]}\fint_{B_\tau}|u|^2dx\leq1,
\end{equation*}
\begin{equation*}
\fint_{Q_1}|\nabla\times \tilde{u}|^2dxdt
=\frac{1}{J^2}\fint_{Q_\tau}|\nabla\times u|^2dxdt\leq1,
\end{equation*}
\begin{equation*}
\begin{aligned}
\fint_{Q_r}|\tilde{f}-\overline{\tilde{f}}_{B_r}|^2dxdt
=\frac{1}{J^2}\fint_{Q_{\tau r}}|f-\overline{f}_{B_{\tau r}}|^2dxdt\leq\varepsilon^2_0r^{2\alpha}
\end{aligned}
\end{equation*}
and
\begin{equation*}
\begin{aligned}
\sup_{t\in(-r^2,0]}\fint_{B_r}|\nabla \tilde{g}|^qdx
=\frac{\tau^q}{J^q}\sup_{t\in(-\tau^2r^2,0]}\fint_{B_{\tau r}}|\nabla g|^qdx\leq1
\end{aligned}
\end{equation*}
for any $0<r\leq1$.
By Theorem 3.1, there exist $\tilde{a}(t)\in L_\infty(-r^2_0,0;\mathbb{R}^{n})$ and $\tilde{b}(t)\in L_\infty(-r^2_0,0;\mathbb{R}^{n^2})$ such that
\begin{equation*}
\fint_{Q_{r}}|\nabla\times \tilde{u}-\overline{(\nabla\times \tilde{u})}_{Q_{r}}|^2dxdt\leq Cr^{2\alpha}
\end{equation*}
and
\begin{equation*}
\begin{aligned}
\sup_{t\in(-r^{2},0]}\fint_{B_{r}}|\tilde{u}-\tilde{a}(t)-\tilde{b}(t)x|^2dx
\leq Cr^{2(1+\alpha)}
\end{aligned}
\end{equation*}
for any $0<r<1$, where $C$ depends only on $n$, $\lambda$, $\alpha$ and $q$.
After rescaling back, we get (\ref{p01}) and (\ref{p02})
with
$a(t)=\tau J\tilde{a}(t)$ and $b(t)=J\tilde{b}(t)$.

\bigskip

We will prove Theorem 3.1 by an iteration staring with the following lemma.

\bigskip

\begin{lemma}
Under the hypotheses of Theorem 3.1,
there exist constants $0<r_0\leq\frac16$, $0<\hat{\varepsilon}_0<1$ and $C_0>0$ depending only on $n$, $\lambda$, $\alpha$ and $q$ such that if
\begin{equation}\label{sh.1}
\big|A(x,t)-A(0,t)\big|\leq\hat{\varepsilon}_0~~\mbox{in}~~Q_1,
\end{equation}
\begin{equation}\label{sh.2}
\sup_{t\in(-1,0]}\fint_{B_1}|u|^2dx\leq1,~~\fint_{Q_1}|\nabla\times u|^2dxdt\leq1,
\end{equation}
\begin{equation}\label{sh.3}
\fint_{Q_1}|f|^2dxdt\leq\hat{\varepsilon}^2_0
~~\mbox{and}~~\sup_{t\in(-1,0]}\fint_{B_1}|\nabla g|^qdx\leq1,
\end{equation}
then there is $b_0(t)\in L_\infty(-r^2_0,0;\mathbb{R}^{n^2})$ satisfying
\begin{equation}\label{sh.21}
\sup_{t\in(-r_0^2,0]}|b_0(t)|\leq C_0,
\end{equation}
\begin{equation}\label{sh.4}
\fint_{Q_{r_0}}|\nabla\times u-\overline{(\nabla\times u)}_{Q_{r_0}}|^2dxdt\leq r^{2\alpha}_0
\end{equation}
and
\begin{equation}\label{f.09}
\begin{aligned}
\sup_{t\in(-r_0^2,0]}\fint_{B_{r_0}}|u-\overline{u}_{B_{r_0}}(t)-b_0(t)x|^2dx
\leq r^{2+2\alpha}_0.
\end{aligned}
\end{equation}
\end{lemma}

\vskip 2mm\noindent{\bf Proof.}\quad
Suppose (\ref{sh.1})-(\ref{sh.3}) hold.  We divide the proof of (\ref{sh.21})-(\ref{f.09}) into 7 steps and the constants $r_0$, $\hat{\varepsilon}_0$ and $C_0$ will be determined in Step 7.

\emph{Step 1. Estimates of $\nabla u$.}~
Using Lemma 2.3 (i) to $u$ in $B_1$ with $p=2$, $\zeta=0$ and $\eta=\overline{g}_{B_1}(t)$, we have,
by $\mbox{div}~u=g$,
\begin{eqnarray*}
\begin{aligned}
\int_{B_{\frac34}}|\nabla u|^{2}dx
\leq&\ C\Big(\int_{B_1}|\nabla\times u|^{2}dx+\int_{B_1}|g-\overline{g}_{B_1}(t)|^{2}dx
+\int_{B_1}|u|^{2}dx\Big)
\end{aligned}
\end{eqnarray*}
for $t\in(-\frac{9}{16},0]$, where $C$ depends only on $n$.
Integrating with respect to $t$ and by Poincar\'{e}'s inequality,
\begin{eqnarray*}
\begin{aligned}
\int_{Q_{\frac34}}|\nabla u|^{2}dxdt
\leq&\ C\Big(\int_{Q_1}|\nabla\times u|^{2}dxdt+\int_{Q_1}|\nabla g(t)|^{2}dxdt+\int_{Q_1}|u|^{2}dxdt\Big).
\end{aligned}
\end{eqnarray*}
In view of  (\ref{sh.2}) and (\ref{sh.3}), we obtain
\begin{eqnarray*}
\begin{aligned}
\int_{Q_{\frac34}}|\nabla u|^{2}dxdt
\leq C
\end{aligned}
\end{eqnarray*}
for some $C$ depending only on $n$.

From (\ref{s0.0}), one has
\begin{eqnarray}\label{sh.53}
\left\{
\begin{array}{ll}
\vspace{2mm}
u_t-\mbox{div}\big(A(0,t)\nabla u\big)+\nabla p=\mbox{div}\Big(\big(A(x,t)
-A(0,t)\big)\nabla u+f\Big)~~\mbox{in}~~Q_1,\\
\vspace{2mm}
\mbox{div}~u=g~~\mbox{in}~~Q_1.
\end{array}
\right.
\end{eqnarray}
Let $h(x,t)=\big(A(x,t)-A(0,t)\big)\nabla u$. From (\ref{sh.1}), it follows that
\begin{eqnarray}\label{sh.7}
\begin{aligned}
\int_{Q_{\frac34}}|h|^{2}dxdt
= \int_{Q_{\frac34}}|\big(A(x,t)-A(0,t)\big)\nabla u|^{2}dxdt
\leq C^*\hat{\varepsilon}_0^2,
\end{aligned}
\end{eqnarray}
where $C^*$ depends only on $n$.

\emph{Step 2. Proof of (\ref{sh.4}).}
Taking $\nabla\times$ on both sides of the first equation in (\ref{sh.53}) and by $\nabla\times\nabla p=0$, we obtain
\begin{equation*}
(\nabla\times u)_t-\mbox{div}\big(A(0,t)\nabla(\nabla\times u)\big)
=\nabla\times\mbox{div}\big(h+f\big)~~\mbox{in}~~Q_1.
\end{equation*}

Decompose $\nabla\times u$ into $w+v$ with $w$ solving
\begin{equation}\label{h.10}
\left\{
\begin{array}{ll}
w_t-\mbox{div}\big(A(0,t)\nabla w\big)
=\nabla\times \mbox{div}~\big(h+f\big)~~\mbox{in}~~Q_{\frac34},\\
\vspace{2mm}
w=0~~\mbox{on}~~\partial_pQ_{\frac34}.
\end{array}
\right.
\end{equation}
Then it follows that
\begin{equation}\label{h.8}
\left\{
\begin{array}{ll}
v_t-\mbox{div}\big(A(0,t)\nabla v\big)
=0~~\mbox{in}~~Q_{\frac34},\\
\vspace{2mm}
v=\nabla\times u~~\mbox{on}~~\partial_pQ_{\frac34}.
\end{array}
\right.
\end{equation}

Applying Lemma 2.4 (ii) to (\ref{h.10}) and in view of (\ref{sh.3}) and (\ref{sh.7}), we have
\begin{eqnarray}\label{sh.9}
\begin{aligned}
\int_{Q_{\frac12}}|w|^2dxdt
\leq C\int_{Q_{\frac34}}|h+f|^2dxdt
\leq C\hat{\varepsilon}_0^2
\end{aligned}
\end{eqnarray}
and then
\begin{eqnarray}\label{sh.10}
\begin{aligned}
\fint_{Q_{r_0}}|w-\overline{w}_{Q_{r_0}}|^2dxdt
\leq\frac{C}{r^{n+2}_0}\int_{Q_{\frac12}}|w|^2dxdt
\leq C_1\frac{\hat{\varepsilon}_0^2}{r^{n+2}_0}
\end{aligned}
\end{eqnarray}
for some $C_1$ depending only on $n$ and $\lambda$.

Apply Lemma 2.4 (i) to (\ref{h.8}) and then
\begin{eqnarray*}
\begin{aligned}
\fint_{Q_{r_0}}|v-\overline{v}_{Q_{r_0}}|^2dxdt
\leq&\ \|v-\overline{v}_{Q_{r_0}}\|^2_{L_\infty(Q_{r_0})}\\
\leq&\ Cr^2_0\Big(\|Dv\|_{L_\infty(Q_{\frac14})}+\|v_t\|_{L_\infty(Q_{\frac14})}\Big)^2\leq Cr^2_0\|v\|^2_{L_2(Q_{\frac12})},
\end{aligned}
\end{eqnarray*}
where $C$ depends only on $n$ and $\lambda$.
From (\ref{sh.2}),  (\ref{sh.9}) and $v=\nabla\times u-w$, it follows that

\begin{equation}\label{jjf}
\begin{aligned}
\fint_{Q_{r_0}}|v-\overline{v}_{Q_{r_0}}|^2dxdt
\leq&\ Cr^2_0\Big(\|\nabla\times u\|^2_{L_2(Q_{\frac12})}+\|w\|^2_{L_2(Q_{\frac12})}\Big)\leq C_2r^2_0,
\end{aligned}
\end{equation}
where $C_2$ depends only on $n$ and $\lambda$.

Combining with (\ref{sh.10}), we get
\begin{equation}\label{sh.8}
\begin{aligned}
\lefteqn{\fint_{Q_{r_0}}|\nabla\times u-\overline{(\nabla\times u)}_{Q_{r_0}}|^2dxdt}\hspace*{18mm}\\
\leq&\ 2\fint_{Q_{r_0}}|w-\overline{w}_{Q_{r_0}}|^2dxdt+2\fint_{Q_{r_0}}|v-\overline{ v}_{Q_{r_0}}|^2dxdt\\
\leq&\ 2C_1\frac{\hat{\varepsilon}_0^2}{r^{n+2}_0}+2C_2r^{2}_0.
\end{aligned}
\end{equation}

\emph{Step 3. Decomposition.} We decompose $u$ and $p$, satisfying (\ref{sh.53}), by the following way. Let $(\hat w, p_{\hat w})$ be the solution of :
\begin{eqnarray}\label{sh.11}
\left\{
\begin{array}{ll}
\vspace{2mm}
\hat{w}_t-\mbox{div}\big(A(0,t)\nabla \hat{w}\big)+\nabla p_{\hat{w}}=\mbox{div}\big(h+f\big)~~\mbox{in}~~Q_{\frac12},\\
\vspace{2mm}
\mbox{div}~\hat{w}=0~~\mbox{in}~~Q_{\frac12},\\
\vspace{2mm}
\hat{w}=0~~\mbox{on}~~\partial_pQ_{\frac12}.
\end{array}
\right.
\end{eqnarray}
Set $\hat v=u-\hat w$ and $p_{\hat v}=p-p_{\hat w}$. It follows that
\begin{eqnarray}\label{sh.12}
\left\{
\begin{array}{ll}
\vspace{2mm}
\hat{v}_t-\mbox{div}\big(A(0,t)\nabla \hat{v}\big)+\nabla p_{\hat{v}}=0~~\mbox{in}~~Q_{\frac12},\\
\vspace{2mm}
\mbox{div}~\hat{v}=g~~\mbox{in}~~Q_{\frac12},\\
\vspace{2mm}
\hat{v}=u~~\mbox{on}~~\partial_pQ_{\frac12}.
\end{array}
\right.
\end{eqnarray}

\emph{Step 4. Estimates of $\hat{w}$.} Multiplying the first equation in (\ref{sh.11}) by $\hat{w}$ and then taking integral with respect to $x$, we deduce, by using the other two equations in (\ref{sh.11}),
\begin{eqnarray*}
\begin{aligned}
\frac{d}{dt}\int_{B_{\frac12}}|\hat{w}|^2dx+\int_{B_{\frac12}}|\nabla\hat{w}|^2dx
\leq C\int_{B_{\frac12}}|h+f|^2dx
\end{aligned}
\end{eqnarray*}
for $t\in(-\frac{1}{4},0]$. Next taking integral with respect to $t$,
\begin{eqnarray}\label{sh.16}
\begin{aligned}
\sup_{t\in(-\frac14,0]}\int_{B_{\frac12}}|\hat{w}|^2dx+\int_{Q_{\frac12}}|\nabla\hat{w}|^2dxdt
\leq C\int_{Q_{\frac12}}|h+f|^2dxdt
\leq C\hat{\varepsilon}_0^2,
\end{aligned}
\end{eqnarray}
where (\ref{sh.3}) and (\ref{sh.7}) are used to derive the last inequality and $C$ depends only on $n$ and $\lambda$.
This implies
\begin{eqnarray}\label{sh.18}
\begin{aligned}
\sup_{t\in(-\frac14,0]}\fint_{B_{r_0}}|\hat{w}-\overline{\hat{w}}_{B_{r_0}}(t)|^2dx
\leq \frac{C}{r^{n}_0}\sup_{t\in(-\frac14,0]}\int_{B_{\frac12}}|\hat{w}|^2dx
\leq C_3\frac{\hat{\varepsilon}_0^2}{r^{n}_0},
\end{aligned}
\end{eqnarray}
where $C_3$ depends only on $n$ and $\lambda$.

\emph{Step 5. Estimates of $\hat{v}$.} It follows from $\hat{v}=u-\hat{w}$, (\ref{sh.2}) and (\ref{sh.16}) that
\begin{eqnarray}\label{sh.17}
\begin{aligned}
\sup_{t\in(-\frac14,0]}\|\hat{v}\|_{L_2(B_{\frac12})}
\leq \sup_{t\in(-\frac14,0]}\|u\|_{L_2(B_{\frac12})}+\sup_{t\in(-\frac14,0]}\|\hat{w}\|_{L_2(B_{\frac12})}\leq C,
\end{aligned}
\end{eqnarray}
where $C$ depends only on $n$ and $\lambda$.

Take $\nabla\times$ on both sides of the first equation in (\ref{sh.12}) and in view of $\nabla\times\nabla p_{\hat v}=0$ , one has
\begin{eqnarray*}
\begin{array}{ll}
\vspace{2mm}
\big(\nabla\times\hat{v}\big)_t-\mbox{div}\Big(A(0,t)\nabla \big(\nabla\times\hat{v}\big)\Big)=0~~\mbox{in}~~Q_{\frac13}.
\end{array}
\end{eqnarray*}
From Lemma 2.4 (i) and $\nabla\times\hat{v}=\nabla\times u-\nabla\times\hat{w}$,  we deduce
\begin{eqnarray*}
\begin{aligned}
&\sup_{t\in(-\frac{1}{16},0]}\|\nabla\times \hat{v}\|_{L_2(B_{\frac14})}+
\sup_{t\in(-\frac{1}{16},0]}\|D(\nabla\times \hat{v})\|_{L_q(B_{\frac14})}\\
&\mbox{}\hskip0.8cm\leq C\left(\|\nabla\times \hat{v}\|_{L_\infty(Q_{\frac14})}+
\|D(\nabla\times \hat{v})\|_{L_\infty(Q_{\frac14})}\right)&\\
&\mbox{}\hskip2cm\leq C\|\nabla\times \hat{v}\|_{L_2(Q_{\frac12})}
\leq C\Big(\|\nabla\times u\|_{L_2(Q_{\frac12})}+\|\nabla \hat{w}\|_{L_2(Q_{\frac12})}\Big).
\end{aligned}
\end{eqnarray*}
From (\ref{sh.2}) and (\ref{sh.16}), it follows that
\begin{eqnarray}\label{sh.91}
\begin{aligned}
\sup_{t\in(-\frac{1}{16},0]}\|\nabla\times \hat{v}\|_{L_2(B_{\frac14})}+
\sup_{t\in(-\frac{1}{16},0]}\|D(\nabla\times \hat{v})\|_{L_q(B_{\frac14})}
\leq C
\end{aligned}
\end{eqnarray}
for some $C$ depending only on $n$ and $\lambda$.

Apply Lemma 2.3 (i) to $\hat{v}$ with $p=2$, $\zeta=0$ and $\eta=\overline{g}_{B_1}(t)$.
By $\mbox{div}~\hat{v}=g$,
we have
\begin{eqnarray}\label{sh.60}
\begin{aligned}
&\sup_{t\in(-\frac{1}{36},0]}\|\nabla\hat{v}\|_{L_2(B_{\frac16})}
\leq C\sup_{t\in(-\frac{1}{16},0]}\Big(\|\nabla\times \hat{v}\|_{L_2(B_{\frac14})}+\|g-\overline{g}_{B_1}(t)\|_{L_2(B_{\frac14})}
+\|\hat{v}\|_{L_2(B_{\frac14})}\Big)\\
&\mbox{}\hskip0.8cm\leq C\sup_{t\in(-\frac{1}{16},0]}\Big(\|\nabla\times \hat{v}\|_{L_2(B_{\frac14})}+\|\nabla g(t)\|_{L_2(B_{\frac14})}
+\|\hat{v}\|_{L_2(B_{\frac14})}\Big)
\leq C,
\end{aligned}
\end{eqnarray}
where the second inequality is derived from Poincar\'{e}'s inequality and the last inequality is derived from
 (\ref{sh.3}), (\ref{sh.17}) and (\ref{sh.91}).

Applying Lemma 2.3 (ii) to $\hat{v}$ with $p=q$ and by (\ref{sh.3}), (\ref{sh.17}) and (\ref{sh.91}), we have
\begin{eqnarray*}
\begin{aligned}
\sup_{t\in(-\frac{1}{36},0]}\|D^2\hat{v}\|_{L_q(B_{\frac16})}
\leq&\ C\sup_{t\in(-\frac{1}{16},0]}\Big(\|D(\nabla\times \hat{v})\|_{L_q(B_{\frac14})}+\|\nabla g\|_{L_q(B_{\frac14})}
+\|\hat{v}\|_{L_2(B_{\frac14})}\Big)
\leq C,
\end{aligned}
\end{eqnarray*}
which leads to
\begin{eqnarray}\label{sh.19}
\begin{aligned}
\lefteqn{\sup_{t\in(-r_0^2,0]}\fint_{B_{r_0}}|\hat{v}-\overline{\hat{v}}_{B_{r_0}}(t)-\overline{\nabla \hat{v}}_{B_{r_0}}(t)x|^2dx}\hspace*{16mm}\\
\leq&\ Cr^4_0\sup_{t\in(-r_0^2,0]}\fint_{B_{r_0}}|D^2\hat{v}|^2dx
\leq Cr^4_0\sup_{t\in(-r_0^2,0]}\Big(\fint_{B_{r_0}}|D^2\hat{v}|^qdx\Big)^{\frac{2}{q}}\\
\leq&\ Cr^{4-\frac{2n}{q}}_0\sup_{t\in(-r_0^2,0]}\|D^2\hat{v}\|^2_{L_q(B_{\frac16})}
\leq C_4r^{4-\frac{2n}{q}}_0,
\end{aligned}
\end{eqnarray}
where $C_4$ depends only on $n$, $\lambda$ and $q$.

\emph{Step 6. Proof of (\ref{sh.21}) and (\ref{f.09}).} Set $b_0(t)=\overline{\nabla \hat{v}}_{B_{r_0}}(t)$.
It follows from (\ref{sh.60}) that
\begin{equation}\label{n.14}
\begin{aligned}
\sup_{t\in(-r_0^2,0]}|b_0(t)|
\leq&\ \sup_{t\in(-r_0^2,0]}\Big(\fint_{B_{r_0}}|\nabla\hat{v}|^2dx\Big)^{\frac12}\\
\leq &\ \frac{C}{r^{\frac n2}_0}\sup_{t\in(-r_0^2,0]}\Big(\int_{B_{\frac16}}|\nabla\hat{v}|^2dx\Big)^{\frac12}
\leq C_5r^{-\frac n2}_0.
\end{aligned}
\end{equation}
where $C_5$ depends only on $n$ and $\lambda$.

Combining (\ref{sh.18}) and (\ref{sh.19}), we get
\begin{eqnarray}\label{sh.54}
\begin{aligned}
&\sup_{t\in(-r_0^2,0]}\fint_{B_{r_0}}|u-\overline{u}_{B_{r_0}}(t)-b_0(t)x|^2dx~~~~~(\mbox{by}~u=\hat v+\hat w)\\
&\leq 2\sup_{t\in(-r_0^2,0]}\left(\fint_{B_{r_0}}|\hat{v}-\overline{\hat{v}}_{B_{r_0}}(t)-\overline{\nabla \hat{v}}_{B_{r_0}}(t)x|^2dx
+\fint_{B_{r_0}}|\hat{w}-\overline{\hat{w}}_{B_{r_0}}(t)|^2dx\right)\\
&\leq 2C_4r^{4-\frac{2n}{q}}_0+2C_3\frac{\hat{\varepsilon}_0^2}{r^{n}_0}.
\end{aligned}
\end{eqnarray}

\emph{Step 7. Determination of the constants.}
Let $C_i$ ($i=1\sim 5$)  be given by
(3.19), (3.20), (3.25), (3.29) and (3.30)
respectively.  Recall $0<\alpha<1$ and $q>\frac{n}{1-\alpha}$. We first set $0<r_0\leq\frac16$ small enough such that
\begin{eqnarray}\label{n6.14}
2C_2r^2_0\leq \frac12r_0^{2\alpha}~~\mbox{and}~~2C_4r^{4-\frac{2n}{q}}_0\leq\frac12r_0^{2+2\alpha}.
\end{eqnarray}
Next, we set $0<\hat{\varepsilon}_0<1$ small enough such that
\begin{eqnarray}\label{n7.14}
2C_1\hat{\varepsilon}_0^2r^{-n-2}_0\leq \frac12r_0^{2\alpha}~~\mbox{and}~~2C_3\hat{\varepsilon}_0^2{r^{-n}_0}\leq\frac12r_0^{2+2\alpha}.
\end{eqnarray}
Then (\ref{sh.4}) and (\ref{f.09}) follow from (\ref{sh.8}) and (\ref{sh.54}) respectively.
Finally, set $C_0=C_5r^{-\frac n2}_0$ and then (\ref{sh.21}) follows from (\ref{n.14}).
\qed

\bigskip

\begin{lemma}
Suppose the hypotheses of Theorem 3.1 hold. Let $r_0$, $\hat{\varepsilon}_0$ and $C_0$ be given by Lemma 3.2 and
\begin{equation}\label{jj2}
\hat{C}_0=C_0+ \bar{C}_0r^{-\frac n2}_0\sum_{k=1}^\infty r^{k\alpha}_0
\end{equation}
with some constant $\bar{C}_0>0$ depending only on $n$ and $\lambda$.
There exists a constant $0<\varepsilon_0\leq\hat{\varepsilon}_0$ depending only on $n$, $\lambda$, $\alpha$ and $q$ such that if (\ref{sh.23})-(\ref{sh.25}) hold, then for any $k\geq1$,
there are constants $C_k\leq\hat{C}_0$ and $b_k(t)\in L_\infty(-r^{2k}_0,0;\mathbb{R}^{n^2})$ satisfying
\begin{equation}\label{sh.27}
\sup_{t\in(-r^{2k}_0,0]}|b_k(t)|\leq C_k,
\end{equation}
\begin{equation}\label{sh.26}
\fint_{Q_{r^k_0}}|\nabla\times u-\overline{(\nabla\times u)}_{Q_{r^k_0}}|^2dxdt\leq r^{2k\alpha}_0
\end{equation}
and
\begin{equation}\label{sh.28}
\begin{aligned}
\sup_{t\in(-r^{2k}_0,0]}\fint_{B_{r^k_0}}|u-\overline{u}_{B_{r^k_0}}(t)-b_k(t)x|^2dx
\leq r^{2k(1+\alpha)}_0.
\end{aligned}
\end{equation}
\end{lemma}

\vskip 2mm\noindent{\bf Proof.}\quad
We prove the lemma by induction on $k$.
Clearly, it holds as $k=1$ with $C_1=C_0$  and $b_1(t)=b_0(t)$ by Lemma 3.2.
Now we suppose there exist $C_k$ and $b_{k}(t)$ such that (\ref{sh.27})-(\ref{sh.28}) hold for some $k\geq1$ and dispose of the case $k+1$ by 7 steps. The constants $\varepsilon_0$ and $\bar C_0$ will be determined in Step 1 and Step 5 respectively.

\emph{Step 0. Normalization.} For $(x,t)\in Q_1$, define
\begin{equation*}
\tilde{u}(x,t)=r_0^{-k(1+\alpha)}u(r^k_0x,r^{2k}_0t),~~\tilde{p}(x,t)=r_0^{-k\alpha}p(r^k_0x,r^{2k}_0t),
\end{equation*}
\begin{equation*}
\tilde{A}(x,t)=A(r^k_0x,r^{2k}_0t),
\end{equation*}
\begin{equation*}
\tilde{f}(x,t)=r_0^{-k\alpha}\big(f(r^k_0x,r^{2k}_0t)-\overline{f}_{B_{r^k_0}}(r^{2k}_0t)\big)~~\mbox{and}~~\tilde{g}(x,t)=r_0^{-k\alpha}g(r^k_0x,r^{2k}_0t).
\end{equation*}
It follows that
\begin{eqnarray}\label{sh.58}
\left\{
\begin{array}{ll}
\vspace{2mm}
\tilde{u}_t-\mbox{div}\big(\tilde{A}(0,t)\nabla \tilde{u}\big)+\nabla \tilde{p}=\mbox{div}\Big(\tilde{f}+\big(\tilde{A}(x,t)-\tilde{A}(0,t)\big)\nabla \tilde{u}\Big)~~\mbox{in}~~Q_{1},\\
\vspace{2mm}
\mbox{div}~\tilde{u}=\tilde{g}~~\mbox{in}~~Q_{1}.
\end{array}
\right.
\end{eqnarray}

By (\ref{sh.23}), (\ref{n.25}) and (\ref{sh.25}),
\begin{equation}\label{sh.35}
\begin{aligned}
|\tilde{A}(x,t)-\tilde{A}(0,t)|=|A(r^k_0x,r^{2k}_0t)-A(0,r^{2k}_0t)|
\leq\varepsilon_0r_0^{k\alpha},
\end{aligned}
\end{equation}
\begin{eqnarray}\label{sh.40}
\begin{aligned}
\fint_{Q_1}|\tilde{f}|^2dxdt
=\fint_{Q_{r^k_0}}|r_0^{-k\alpha}\big(f(x,t)-\overline{f}_{B_{r^k_0}}(t)\big)|^2dxdt
\leq\varepsilon^2_0\leq\hat{\varepsilon}^2_0
\end{aligned}
\end{eqnarray}
and
\begin{equation}\label{sh.90}
\sup_{t\in(-1,0]}\fint_{B_1}|\nabla \tilde{g}|^qdx=r^{kq(1-\alpha)}_0\sup_{t\in(-r^{2k}_0,0]}\fint_{B_{r^k_0}}|\nabla g(x,t)|^qdx\leq1.
\end{equation}

In view of the induction hypotheses (\ref{sh.26}) and (\ref{sh.28}), we have
\begin{equation}\label{sh.30}
\fint_{Q_1}|\nabla\times \tilde{u}-\overline{(\nabla\times \tilde{u})}_{Q_1}|^2dxdt\leq1
\end{equation}
and
\begin{equation}\label{sh.31}
\sup_{t\in(-1,0]}\fint_{B_1}|\tilde{u}-\alpha(t)-\beta(t)x|^2dx\leq1,
\end{equation}
where
\begin{equation*}
\alpha(t)=r^{-k(1+\alpha)}_0\overline{u}_{B_{r^k_0}}(r^{2k}_0t)~~\mbox{and}~~\beta(t)=r^{-k\alpha}_0b_k(r^{2k}_0t).
\end{equation*}

\emph{Step 1. Estimates of $\nabla\tilde u$.}~~
Applying Lemma 2.3 (i) to $\tilde{u}-\alpha(t)-\beta(t)x$ with $p=2$, $\zeta=\overline{(\nabla\times \tilde{u})}_{Q_{1}}-\nabla\times\big(\beta(t)x\big)$ and $\eta=\overline{\tilde{g}}_{B_{1}}(t)-div\big(\beta(t)x\big)$,
we obtain, by $\mbox{div}~\tilde{u}=\tilde{g}$,
\begin{eqnarray*}
\begin{aligned}
\lefteqn{\int_{B_{\frac34}}|\nabla\tilde{u}-\beta(t)|^{2}dx= \int_{B_{\frac34}}|\nabla\big(\tilde{u}-\alpha(t)-\beta(t)x\big)|^{2}dx}\hspace*{3mm}\\
\leq&\ C\Big(\int_{B_1}|\nabla\times \tilde{u}-\overline{(\nabla\times \tilde{u})}_{Q_{1}}|^{2}dx+\int_{B_1}|\tilde{g}-\overline{\tilde{g}}_{B_1}(t)|^{2}dx
+\int_{B_1}|\tilde{u}-\alpha(t)-\beta(t)x|^{2}dx\Big)
\end{aligned}
\end{eqnarray*}
for $t\in(-\frac{9}{16},0]$ and then
integrating it with respect to $t$ and using Poincar\'{e}'s inequality,
\begin{eqnarray*}
\begin{aligned}
&\int_{Q_{\frac34}}|\nabla\tilde{u}-\beta(t)|^{2}dxdt\\
&\leq C\Big(\int_{Q_1}|\nabla\times \tilde{u}-\overline{(\nabla\times \tilde{u})}_{Q_{1}}|^{2}dxdt+\int_{Q_1}|\nabla\tilde{g}|^{2}dxdt
+\int_{Q_1}|\tilde{u}-\alpha(t)-\beta(t)x|^{2}dxdt\Big),
\end{aligned}
\end{eqnarray*}
where $C$ depends only on $n$.
From (\ref{sh.90}), (\ref{sh.30}) and (\ref{sh.31}),
it follows that
\begin{eqnarray*}
\begin{aligned}
\int_{Q_{\frac34}}|\nabla\tilde{u}-\beta(t)|^{2}dxdt
\leq C,
\end{aligned}
\end{eqnarray*}
where $C$ depends only on $n$.
It and (\ref{sh.27}) lead to
\begin{eqnarray}\label{j002}
\begin{aligned}
\int_{Q_{\frac34}}|\nabla\tilde{u}|^{2}dxdt
\leq&\ 2\int_{Q_{\frac34}}|\beta(t)|^{2}dxdt+2\int_{Q_{\frac34}}|\nabla\tilde{u}-\beta(t)|^{2}dxdt\\
&\mbox{}\hskip-1.8cm\leq C\sup_{t\in(-\frac{9}{16},0]}|r^{-k\alpha}_0b_k(r^{2k}_0t)|^{2}+C
\leq Cr^{-2k\alpha}_0\sup_{t\in(-r^{2k}_0,0]}|b_k(s)|^{2}+C\\
\leq&\ CC_kr^{-2k\alpha}_0+C\leq C\big(\hat{C}_0r^{-2k\alpha}_0+1\big).
\end{aligned}
\end{eqnarray}

Set $\tilde{h}(x,t)=\big(\tilde{A}(x,t)-\tilde{A}(0,t)\big)\nabla \tilde{u}$ and then (\ref{sh.58}) reads as
\begin{eqnarray}\label{j003}
\left\{
\begin{array}{ll}
\vspace{2mm}
\tilde{u}_t-\mbox{div}\big(\tilde{A}(0,t)\nabla \tilde{u}\big)+\nabla \tilde{p}=\mbox{div}\big(\tilde{f}+\tilde h\big)~~\mbox{in}~~Q_{1},\\
\vspace{2mm}
\mbox{div}~\tilde{u}=\tilde{g}~~\mbox{in}~~Q_{1}.
\end{array}
\right.
\end{eqnarray}
Combining (\ref{sh.35}) and (\ref{j002}), we get
\begin{eqnarray}\label{sh.36}
\begin{aligned}
\int_{Q_{\frac34}}|\tilde{h}|^{2}dxdt
=&\ \int_{Q_{\frac34}}|\big(\tilde{A}(x,t)-\tilde{A}(0,t)\big)\nabla \tilde{u}|^{2}dxdt\\
\leq &C\varepsilon_0^2r^{2k\alpha}_0\big(\hat{C}_0r^{-2k\alpha}_0+1\big)
\leq C\varepsilon_0^2\big(\hat{C}_0+1\big)\leq C^*\hat{\varepsilon}_0^2
\end{aligned}
\end{eqnarray}
as $0<\varepsilon_0\leq\hat{\varepsilon}_0$ is sufficiently small,
where the constant $C^*$ is the same as in (\ref{sh.7}).

\emph{Step 2. Proof of (\ref{sh.26}).}
Taking $\nabla\times$ on both sides of the first equation in (\ref{j003}) and by the sane arguments as in Step 2 of the proof of Lemma 3.2 with (\ref{sh.40}), (\ref{sh.30}) and (\ref{sh.36}) replacing (\ref{sh.3}), (\ref{sh.2}) and (\ref{sh.7}) respectively, we obtain that
\begin{eqnarray*}
\fint_{Q_{r_0}}|\nabla\times \tilde{u}-\overline{(\nabla\times \tilde{u})}_{Q_{r_0}}|^2dxdt
\leq 2C_1\frac{\hat{\varepsilon}_0^2}{r^{n+2}_0}+2C_2r^{2}_0\leq r_0^{2\alpha},
\end{eqnarray*}
where $C_1$ and $C_2$ are the same as in (\ref{sh.8}) and the last inequality is derived from (\ref{n6.14}) and (\ref{n7.14}).
After rescaling back,
\begin{eqnarray*}
\begin{aligned}
\fint_{Q_{r^{k+1}_0}}|\nabla\times u-\overline{(\nabla\times u)}_{Q_{r^{k+1}_0}}|^2dxdt\leq r_0^{2(k+1)\alpha}.
\end{aligned}
\end{eqnarray*}
That is, (\ref{sh.26}) is proved for $k+1$.

\emph{Step 3. Decomposition.} Similar to Step 3 in the proof of Lemma 3.2, we decompose solution of (\ref{j003}) such that
$(\tilde{u},\tilde p)=(\hat{w},\hat p_{\hat w})+(\hat{v},\hat p_{\hat v})$
with
\begin{eqnarray}\label{sh.38}
\left\{
\begin{array}{ll}
\vspace{2mm}
\hat{w}_t-\mbox{div}\big(\tilde{A}(0,t)\nabla \hat{w}\big)+\nabla p_{\hat{w}}=\mbox{div}\big(\tilde{h}+\tilde{f}\big)~~\mbox{in}~~Q_{\frac12},\\
\vspace{2mm}
\mbox{div}~\hat{w}=0~~\mbox{in}~~Q_{\frac12},\\
\vspace{2mm}
\hat{w}=0~~\mbox{on}~~\partial_pQ_{\frac12}
\end{array}
\right.
\end{eqnarray}
and
\begin{eqnarray}\label{sh.39}
\left\{
\begin{array}{ll}
\vspace{2mm}
\hat{v}_t-\mbox{div}\big(\tilde{A}(0,t)\nabla \hat{v}\big)+\nabla p_{\hat{v}}=0~~\mbox{in}~~Q_{\frac12},\\
\vspace{2mm}
\mbox{div}~\hat{v}=\tilde{g}~~\mbox{in}~~Q_{\frac12},\\
\vspace{2mm}
\hat{v}=\tilde{u}~~\mbox{on}~~\partial_pQ_{\frac12}.
\end{array}
\right.
\end{eqnarray}

\emph{Step 4. Estimates of $\hat{w}$.} Multiplying the first equation in (\ref{sh.38}) by $\hat{w}$,
by the same arguments as in Step 4 of Lemma 3.2 with (\ref{sh.40}) and (\ref{sh.36}) replacing
 (\ref{sh.3}) and (\ref{sh.7}) respectively,
we have
\begin{eqnarray}\label{sh.45}
\begin{aligned}
\sup_{t\in(-\frac14,0]}\int_{B_{\frac12}}|\hat{w}|^2dx+\int_{Q_{\frac12}}|\nabla\hat{w}|^2dxdt
\leq C\int_{Q_{\frac12}}|\tilde{h}+\tilde{f}|^2dxdt
\leq C\hat{\varepsilon}_0^2
\end{aligned}
\end{eqnarray}
and
\begin{eqnarray}\label{sh.46}
\begin{aligned}
\sup_{t\in(-\frac14,0]}\fint_{B_{r_0}}|\hat{w}-\overline{\hat{w}}_{B_{r_0}}(t)|^2dx
\leq \frac{C}{r^{n}_0}\sup_{t\in(-\frac14,0]}\int_{B_{\frac12}}|\hat{w}|^2dx
\leq C_3\frac{\hat{\varepsilon}_0^2}{r^{n}_0},
\end{aligned}
\end{eqnarray}
where $C$ depends only on $n$ and $\lambda$, and $C_3$ is given by (\ref{sh.18}).

\emph{Step 5. Estimates of $\hat{v}$.} By (\ref{sh.31}), (\ref{sh.45}) and $\hat{v}=\tilde{u}-\hat{w}$, it follows that
\begin{eqnarray}\label{sh.47}
\begin{aligned}
 \lefteqn{\sup_{t\in(-\frac14,0]}\|\hat{v}-\alpha(t)-\beta(t)x\|_{L_2(B_{\frac12})}}\hspace*{14mm}\\
\leq&\ \sup_{t\in(-\frac14,0]}\|\tilde{u}-\alpha(t)-\beta(t)x\|_{L_2(B_{\frac12})}+\sup_{t\in(-\frac14,0]}\|\hat{w}\|_{L_2(B_{\frac12})}\leq C,
\end{aligned}
\end{eqnarray}
where $C$ depends only on $n$ and $\lambda$.

Take $\nabla\times$ on both sides of the first equation in (\ref{sh.39}). From
$\nabla\times\nabla p_{\hat v}=0$, we deduce
\begin{eqnarray*}
\begin{array}{ll}
\vspace{2mm}
\big(\nabla\times\hat{v}\big)_t-\mbox{div}\Big(\tilde{A}(0,t)\nabla \big(\nabla\times\hat{v}\big)\Big)=0~~\mbox{in}~~Q_{\frac13}
\end{array}
\end{eqnarray*}
or, observing that $\overline{(\nabla\times \tilde{u})}_{Q_1}$ is a constant,
\begin{eqnarray*}
\begin{array}{ll}
\vspace{2mm}
\big(\nabla\times\hat{v}-\overline{(\nabla\times \tilde{u})}_{Q_1}\big)_t-\mbox{div}\Big(\tilde{A}(0,t)\nabla \big(\nabla\times\hat{v}-\overline{(\nabla\times \tilde{u})}_{Q_1}\big)\Big)=0~~\mbox{in}~~Q_{\frac13}.
\end{array}
\end{eqnarray*}
By Lemma 2.4 (i), (\ref{sh.30}), (\ref{sh.45}) and
$\nabla\times\hat{v}=\nabla\times \tilde{u}-\nabla\times\hat{w}$,
\begin{eqnarray}\label{sh.61}
\begin{aligned}
&\mbox{}\hskip-0.8cm\sup_{t\in(-\frac{1}{16},0]}\|\nabla\times \hat{v}-\overline{(\nabla\times \tilde{u})}_{Q_1}\|_{L_2(B_{\frac14})}+
\sup_{t\in(-\frac{1}{16},0]}\|D(\nabla\times \hat{v})\|_{L_q(B_{\frac14})}\\
&\leq C\left(\|\nabla\times \hat{v}-\overline{(\nabla\times \tilde{u})}_{Q_1}\|_{L_\infty(Q_{\frac14})}
+\|D(\nabla\times \hat{v})\|_{L_\infty(Q_{\frac14})}\right)\\
&\leq C\|\nabla\times \hat{v}-\overline{(\nabla\times \tilde{u})}_{Q_1}\|_{L_2(Q_{\frac12})}
=
C\|\nabla\times \tilde{u}-\overline{(\nabla\times \tilde{u})}_{Q_1}-\nabla\times\hat{w}\|_{L_2(Q_{\frac12})}\\
&\leq C\Big(\|\nabla\times \tilde{u}-\overline{(\nabla\times \tilde{u})}_{Q_1}\|_{L_2(Q_{\frac12})}+\|\nabla\hat{w}\|_{L_2(Q_{\frac12})}\Big)\leq C,
\end{aligned}
\end{eqnarray}
where $C$ depends only on $n$ and $\lambda$.

Apply Lemma 2.3 (i) to $\hat{v}-\alpha(t)-\beta(t)x$ with $p=2$, $\zeta=\overline{(\nabla\times \tilde{u})}_{Q_{1}}-\nabla\times\big(\beta(t)x\big)$ and $\eta=\overline{\tilde{g}}_{B_{1}}(t)-div\big(\beta(t)x\big)$.
Then we have, by $\mbox{div}~\hat{v}=\tilde{g}$ and Poincar\'{e}'s inequality,
\begin{eqnarray*}
\begin{aligned}
&\sup_{t\in(-\frac{1}{36},0]}\|\nabla\hat{v}-\beta(t)\|_{L_2(B_{\frac16})}\\
&\leq C\sup_{t\in(-\frac{1}{16},0]}\Big(\|\nabla\times \hat{v}-\overline{(\nabla\times \tilde{u})}_{Q_1}\|_{L_2(B_{\frac14})}
+\|\tilde{g}-\overline{\tilde{g}}_{B_{1}}(t)\|_{L_2(B_{\frac14})}+\|\hat{v}-\alpha(t)-\beta(t)x\|_{L_2(B_{\frac14})}\Big)\\
&\leq C\sup_{t\in(-\frac{1}{16},0]}\Big(\|\nabla\times \hat{v}-\overline{(\nabla\times \tilde{u})}_{Q_1}\|_{L_2(B_{\frac14})}
+\|\nabla\tilde{g}(t)\|_{L_2(B_{\frac14})}+\|\hat{v}-\alpha(t)-\beta(t)x\|_{L_2(B_{\frac14})}\Big).
\end{aligned}
\end{eqnarray*}
In view of (\ref{sh.90}), (\ref{sh.47}) and (\ref{sh.61}),
\begin{eqnarray*}
\begin{aligned}
\sup_{t\in(-\frac{1}{36},0]}\|\nabla\hat{v}-\beta(t)\|_{L_2(B_{\frac16})}
\leq C,
\end{aligned}
\end{eqnarray*}
which combined with (\ref{sh.27}) leads to
\begin{eqnarray}\label{sh.62}
\begin{aligned}
\lefteqn{\sup_{t\in(-r_0^2,0]}\Big(\fint_{B_{r_0}}|\nabla\hat{v}|^2dx\Big)^{\frac12}}\hspace*{10mm}\\
\leq&\ \sup_{t\in(-r_0^2,0]}\Big(\fint_{B_{r_0}}|\beta(t)|^2dx\Big)^{\frac12}
+\sup_{t\in(-r_0^2,0]}\Big(\fint_{B_{r_0}}|\nabla\hat{v}-\beta(t)|^2dx\Big)^{\frac12}\\
\leq&\ \sup_{t\in(-r_0^2,0]}|r^{-k\alpha}_0b_k(r^{2k}_0t)|
+\frac{C}{r^{\frac n2}_0}\sup_{t\in(-r_0^2,0]}\Big(\int_{B_{\frac16}}|\nabla\hat{v}-\beta(t)|^2dx\Big)^{\frac12}\\
\leq&\ r^{-k\alpha}_0\sup_{t\in(-r_0^{2+2k},0]}|b_k(s)|+\bar{C}_0r^{-\frac n2}_0
\leq C_kr^{-k\alpha}_0+\bar{C}_0r^{-\frac n2}_0,
\end{aligned}
\end{eqnarray}
where $C$ and $\bar{C}_0$ depend only on $n$ and $\lambda$.

Apply Lemma 2.3 (ii) to $\hat{v}-\alpha(t)-\beta(t)x$ and then by (\ref{sh.90}), (\ref{sh.47}) and (\ref{sh.61}),
\begin{eqnarray*}
\begin{aligned}
\sup_{t\in(-\frac{1}{36},0]}\|D^2\hat{v}\|_{L_q(B_{\frac16})}
\leq&\ C\sup_{t\in(-\frac{1}{16},0]}\Big(\|D\nabla\times \hat{v}\|_{L_q(B_{\frac14})}+\|\nabla \tilde{g}\|_{L_q(B_{\frac14})}\\
&\mbox{}\hskip3cm  +\|\hat{v}-\alpha(t)-\beta(t)x\|_{L_2(B_{\frac14})}\Big)
\leq C.
\end{aligned}
\end{eqnarray*}
It follows that
\begin{eqnarray}\label{sh.50}
\begin{aligned}
\lefteqn{\sup_{t\in(-r_0^2,0]}\fint_{B_{r_0}}|\hat{v}-\overline{\hat{v}}_{B_{r_0}}(t)-\overline{\nabla \hat{v}}_{B_{r_0}}(t)x|^2dx}\hspace*{16mm}\\
\leq&\ Cr^4_0\sup_{t\in(-r_0^2,0]}\fint_{B_{r_0}}|D^2\hat{v}|^2dx
\leq Cr^4_0\sup_{t\in(-r_0^2,0]}\Big(\fint_{B_{r_0}}|D^2\hat{v}|^qdx\Big)^{\frac{2}{q}}\\
\leq&\ Cr^{4-\frac{2n}{q}}_0\sup_{t\in(-r_0^2,0]}\|D^2\hat{v}\|^2_{L_q(B_{\frac16})}
\leq C_4r^{4-\frac{2n}{q}}_0,
\end{aligned}
\end{eqnarray}
where $C_4$ is the same constant as in (\ref{sh.19}).

\emph{Step 6. Proof of (\ref{sh.27}) and (\ref{sh.28}).}
By (\ref{sh.46}) and (\ref{sh.50}), we obtain
\begin{eqnarray*}
\begin{aligned}
\lefteqn{\sup_{t\in(-r_0^2,0]}\fint_{B_{r_0}}|\tilde{u}-\overline{\tilde{u}}_{B_{r_0}}(t)-\overline{\nabla \hat{v}}_{B_{r_0}}(t)x|^2dx}\hspace*{12mm}\\
\leq&\ 2\sup_{t\in(-r_0^2,0]}\fint_{B_{r_0}}|\hat{v}-\overline{\hat{v}}_{B_{r_0}}(t)-\overline{\nabla \hat{v}}_{B_{r_0}}(t)x|^2dx
+2\sup_{t\in(-r_0^2,0]}\fint_{B_{r_0}}|\hat{w}-\overline{\hat{w}}_{B_{r_0}}(t)|^2dx\\
\leq&\ 2C_4r^{4-\frac{2n}{q}}_0+2C_3\frac{\hat{\varepsilon}_0^2}{r^{n}_0}\leq r_0^{2+2\alpha},
\end{aligned}
\end{eqnarray*}
where the last inequality is derived from (\ref{n6.14}) and (\ref{n7.14}).
Rewriting the above inequality,
\begin{equation*}
\begin{aligned}
\sup_{t\in(-r_0^{2(k+1)},0]}\fint_{B_{r^{k+1}_0}}|u-\overline{u}_{B_{r^{k+1}_0}}(t)-b_{k+1}(t)x|^2dx
\leq r^{2(k+1)(1+\alpha)}_0
\end{aligned}
\end{equation*}
with
\begin{eqnarray*}
b_{k+1}(t)=r^{k\alpha}_0\overline{\nabla\hat{v}}_{B_{r_0}}\Big(\frac{t}{r^{2k}_0}\Big).
\end{eqnarray*}
By (\ref{sh.62}), we have
\begin{eqnarray*}
\begin{aligned}
\sup_{t\in(-r_0^{2(k+1)},0]}|b_{k+1}(t)|
\leq&\ \sup_{t\in(-r_0^{2(k+1)},0]}r_0^{k\alpha}\Big(\fint_{B_{r_0}}|\nabla \hat{v}(x,\frac{t}{r^{2k}_0})|^2dx\Big)^{\frac12}\\
=&\ \sup_{t\in(-r_0^{2},0]}r_0^{k\alpha}\Big(\fint_{B_{r_0}}|\nabla \hat{v}|^2dx\Big)^{\frac12}
\leq r_0^{k\alpha}\Big(C_kr^{-k\alpha}_0+\frac{\bar{C}_0}{r^{\frac n2}_0}\Big)\\
\leq&\ C_k+\bar{C}_0r^{-\frac n2}_0r_0^{k\alpha}:=C_{k+1}.
\end{aligned}
\end{eqnarray*}
Then from induction, we deduce
$$C_{k+1}=C_0+{\bar C}_0r_0^{-\frac{n}2}\sum_{j=1}^{k}r_0^{j\alpha}.$$
In view of (\ref{jj2}), we have $C_{k+1}\leq\hat C_0$.
That is, (\ref{sh.27}) and (\ref{sh.28}) are obtained for $k+1$ and the proof is completed.
\qed

\bigskip

\noindent{\bf Proof of Theorem 3.1.}\quad
For any $0<r<1$, let $k\geq1$ such that
\begin{equation}\label{j05}
r^{k+1}_0\leq r\leq r^k_0.
\end{equation}
By (\ref{sh.26}), we have
\begin{equation*}
\begin{aligned}
&\fint_{Q_{r}}|\nabla\times u-\overline{(\nabla\times u)}_{Q_{r}}|^2dxdt\\
&\mbox{}\hskip1cm\leq2\frac{|Q_{r^k_0}|}{|Q_{r}|}\fint_{Q_{r^k_0}}|\nabla\times u-\overline{(\nabla\times u)}_{Q_{r^k_0}}|^2dxdt
\leq Cr^{2k\alpha}_0
\leq Cr^{2\alpha}
\end{aligned}
\end{equation*}
and hence (\ref{sh.76}) holds.

For any $t\in(-r^2_0,0)$, define
\begin{equation*}
a(t)=\overline{u}_{B_{r^k_0}}(t)~\mbox{as} ~t\in(-r^{2k}_0,-r^{2(k+1)}_0]
\end{equation*}
and
\begin{equation*}
b(t)=b_k(t)~\mbox{as}~t\in(-r^{2k}_0,-r^{2(k+1)}_0]
\end{equation*}
for $k\geq 1$, where $b_k$ is given by Lemma 3.3.

From (\ref{sh.28}), we deduce
\begin{equation*}
\begin{aligned}
&\mbox{}\hskip-1.6cm\sup_{t\in(-r^{2(k+1)}_0,0]}\fint_{B_{r_0^{k+1}}}|\overline{u}_{B_{r^{k+1}_0}}(t)-\overline{u}_{B_{r^k_0}}(t)+b_{k+1}(t)x-b_{k}(t)x|^2dx\\
&\leq C\Bigg\{\sup_{t\in(-r^{2k}_0,0]}\fint_{B_{r_0^{k}}}|u-\overline{u}_{B_{r_0^{k}}}(t)-b_{k}(t)x|^2dx\\
&\mbox{}\hskip2cm +\sup_{t\in(-r^{2(k+1)}_0,0]}\fint_{B_{r_0^{k+1}}}|u-\overline{u}_{B_{r_0^{k+1}}}(t)-b_{k+1}(t)x|^2dx\Bigg\}\\
&\leq Cr^{2k(1+\alpha)}_0
\end{aligned}
\end{equation*}
and then by Lemma 2.2, we have
\begin{equation*}
\begin{aligned}
\sup_{t\in(-r^{2(k+1)}_0,0]}|\overline{u}_{B_{r^{k}_0}}(t)-\overline{u}_{B_{r^{k+1}_0}}(t)|
\leq&\ Cr^{k(1+\alpha)}_0
\end{aligned}
\end{equation*}
and
\begin{equation*}
\sup_{t\in(-r^{2(k+1)}_0,0]}|b_{k}(t)-b_{k+1}(t)|\leq Cr^{k\alpha}_0.
\end{equation*}
Therefore we can and we do define
$$a(0)=\lim_{k\rightarrow\infty}\overline{u}_{B_{r^{k}_0}}(0)~~\mbox{and}~~b(0)=\lim_{k\rightarrow\infty}b_{k}(0).$$

Since
\begin{equation*}
\begin{aligned}
\lefteqn{\sup_{t\in(-r^{2k}_0,0]}|a(t)-\overline{u}_{B_{r^k_0}}(t)|
=\sup_{t\in(-r^{2(k+1)}_0,0]}|a(t)-\overline{u}_{B_{r^k_0}}(t)|}\hspace*{10mm}\\
\leq&\ \sup_{t\in(-r^{2(k+1)}_0,0]}\Big(|a(t)-\overline{u}_{B_{r^{k+1}_0}}(t)|+
|\overline{u}_{B_{r^k_0}}(t)-\overline{u}_{B_{r^{k+1}_0}}(t)|
\Big)\\
\leq&\
\sup_{t\in(-r^{2(k+1)}_0,0]}|a(t)-\overline{u}_{B_{r^{k+1}_0}}(t)|+
Cr^{k(1+\alpha)}_0,
\end{aligned}
\end{equation*}
we have, by repeating the above inequalities,
\begin{equation*}
\sup_{t\in(-r^{2k}_0,0]}|a(t)-\overline{u}_{B_{r^k_0}}(t)|\leq C\sum_{j=k}^\infty r_0^{j(1+\alpha)}\leq Cr_0^{k(1+\alpha)}
\end{equation*}
and similarly,
\begin{equation*}
\sup_{t\in(-r^{2k}_0,0]}|b(t)-b_{k}(t)|\leq C\sum_{j=k}^\infty r_0^{j\alpha}\leq Cr_0^{k\alpha},
\end{equation*}
where $C$ depends only on $n$, $\lambda$, $\alpha$ and $q$.
By (\ref{sh.28}) again,
\begin{equation*}
\begin{aligned}
\lefteqn{\sup_{t\in(-r^{2k}_0,0]}\fint_{B_{r_0^k}}|u-a(t)-b(t)x|^2dx}\hspace*{4mm}\\
\leq&\ 3\sup_{t\in(-r^{2k}_0,0]}\Big(\fint_{B_{r^k_0}}|u-\overline{u}_{B_{r^k_0}}(t)-b_k(t)x|^2dx+|a(t)-\overline{u}_{B_{r^k_0}}(t)|^2
+|b(t)-b_{k}(t)|^2r_0^{2k}\Big)\\
\leq&\ Cr_0^{2k(1+\alpha)}.
\end{aligned}
\end{equation*}
For any $0<r\leq r_0$, let $k$ be given by (\ref{j05}) and it follows that
\begin{equation*}
\begin{aligned}
\sup_{t\in(-r^{2},0]}\fint_{B_{r}}|u-a(t)-b(t)x|^2dx
\leq&\ r_0^{-n}\sup_{t\in(-r^{2k}_0,0]}\fint_{B_{r_0^k}}|u-a(t)-b(t)x|^2dx\\
\leq&\ Cr_0^{2k(1+\alpha)}\leq Cr^{2(1+\alpha)}.
\end{aligned}
\end{equation*}
Thus (\ref{f4.09}) is proved.

Finally, we show the boundedness of $a(t)$ and $b(t)$. Actually,
by (\ref{sh.27}), we obtain
\begin{equation*}
\sup_{t\in(-r^{2}_0,0]}|b(t)|\leq\sup_k\sup_{t\in(-r^{2k}_0,-r^{2(k+1)}_0]}|b_{k}(t)|
\leq\sup_k C_k
\leq \hat{C}_0
\end{equation*}
and from (\ref{sh.24}), we see
\begin{equation*}
\begin{aligned}
\sup_{t\in(-r^{2}_0,0]}|a(t)|=&\
\sup_{t\in(-r^{2}_0,0]}|\overline{u}_{B_{r_0}}(t)|+
\sup_{t\in(-r^{2}_0,0]}|a(t)-\overline{u}_{B_{r_0}}(t)|
\\
\leq&\ Cr^{-\frac n2}_0\sup_{t\in(-r^{2}_0,0]}\Big(\int_{B_{1}}|u|^2dx\Big)^{\frac12}+Cr_0^{1+\alpha}.
\end{aligned}
\end{equation*}
The proof of Theorem 3.1 is complete.
\qed

\bigskip

\noindent{\bf Proof of Corollary 1.2.}\quad Let $R=\min\{\mbox{dist}(\Omega',\partial\Omega),\sqrt{T-T'}\}$ and $q=\frac{2n}{1-\alpha}$.
For any $(x_0,t_0)\in\Omega'_{T'}$, we have (\ref{t1}), (\ref{t2}) and (\ref{t3}) hold with $M_1=\|f\|^2_{C_{x}^{\alpha}(\Omega_{T})}$ and $M_2=\|\nabla g\|_{L_\infty(\Omega_T)}^q$. Then by Theorem 1.1, there exists
$\sigma$ depending only on $n$, $\lambda$, $\alpha$ and $M_0$ such that (\ref{p01}) holds for any $0<r\leq r_0$ with $r_0=\sigma R$. From Lemma 2.1 (i), we obtain the estimate of $[\nabla\times u]_{C^{\alpha,\frac{\alpha}{2}}(\overline{\Omega'_{T'}})}$. Then the desired estimate
for $\|\nabla\times u\|_{C^{\alpha,\frac{\alpha}{2}}(\overline{\Omega'_{T'}})}$ in (\ref{g01}) follows from the boundedness of $\|\nabla\times u\|_{L_2(\Omega_{T})}$ immediately.

Now we turn to the estimate for $\|\nabla u\|_{C_x^{\alpha}(\overline{\Omega'_{T'}})}$.
By Theorem 1.1, there exist
$a(x_0,t)\in L_\infty(t_0-r_0^2,t_0;\mathbb{R}^{n})$ and $b(x_0,t)\in L_\infty(t_0-r_0^2,t_0;\mathbb{R}^{n^2})$
such that (\ref{p02}) holds for any $0<r\leq r_0$ with $a(x_0,t)$ and $b(x_0,t)$ replacing $a(t)$ and $b(t)$ respectively. Set $t=t_0$ in the obtained (\ref{p02}) and then
\begin{equation*}
\begin{aligned}
\lefteqn{
\left(\fint_{B_{r}(x_0)}|u(x,t_0)-a(x_0,t_0)-b(x_0,t_0)x|^2dx\right)^{\frac12}}
\hspace*{8mm}\\
\leq&\ Cr^{1+\alpha}\Big(
\|u\|_{L_\infty(-T,0;L_2(\Omega))}+\|\nabla\times u\|_{L_2(\Omega_{T})}+\|f\|_{C_{x}^{\alpha}(\Omega_{T})}
+\|\nabla g\|_{L_\infty(\Omega_T)}\Big)
\end{aligned}
\end{equation*}
for any $0<r\leq r_0$, where $C$ depends only on $n$, $\lambda$, $\alpha$ and $M_0$.
Then by Lemma 2.1 (ii) and the boundedness of $\|u\|_{L_\infty(-T,0;L_2(\Omega))}$, we obtain the desired estimate of $\|\nabla u\|_{C_x^{\alpha}(\overline{\Omega'_{T'}})}$ clearly. This gives (\ref{g01}) completely.
\qed

\bigskip

\bibliographystyle{plain}

\end{document}